\documentclass[a4paper, 10pt]{article}
\usepackage{amssymb}
\usepackage{amstext}
\usepackage{amsmath}
\usepackage{amscd}
\usepackage{latexsym}
\usepackage{theorem}

\theoremheaderfont{\scshape}
\theorembodyfont{\itshape}
\newtheorem{thm}{\bfseries Theorem}[section]
\newtheorem{prop}[thm]{\bfseries Proposition}
\newtheorem{lemma}[thm]{\bfseries Lemma}
\newtheorem{cor}[thm]{\bfseries Corollary}
\theorembodyfont{\rmfamily}
\newtheorem{defn}[thm]{\bfseries Definition}
\newtheorem{ex}[thm]{\bfseries Example}
\newtheorem{rem}[thm]{\bfseries Remark}

\newtheorem{pf}{Proof.}

\numberwithin{equation}{section}

\def\rmod{{R\textrm{-mod}}}
\def\mods{\textrm{mod-}S}
\def\modr{\textrm{mod-}R}
\def\DR{{\mathfrak D}^b (\rmod )}
\def\DS{{\mathfrak D}^b (\mods)}

\def\Hom{\mathrm{Hom}}
\def\RHom{\mathrm{{\bf R}Hom}}
\def\LT{\overset{\bf L}{\otimes }}
\def\Ext{\mathrm{Ext}}

\def\End{\mathrm{End}}
\def\coker{\mathrm{Coker}}

\def\m{\mathfrak m}

\def\g{\mathrm{G}\textrm{-}\mathrm{dim}}
\def\cm{\mathrm{CM}\textrm{-}\mathrm{dim}}

\def\depth{\mathrm{depth\,}}

\def\pd{\mathrm{pd}}

\def\add{\mathrm{add}}

\def\AC{{\cal R}(C)}
\def\rAC{{\cal R}_R(C)}
\def\AsC{{\cal R}_S(C)}
\def\rAR{{\cal R}_R(R)}
\def\rACd{\mathrm{G}_C\textrm{-}\mathrm{dim}_R\ }
\def\AsCd{\mathrm{G}_C\textrm{-}\mathrm{dim}_S\ }
\def\ACd{\mathrm{G}_C\textrm{-}\mathrm{dim}_R\ }
\def\gc{\mathrm{G}_C}
\def\gcc{\mathrm{G}_{\C}}
\def\rACC{{\cal R}_R(\C)}
\def\AsCC{{\cal R}_S(\C)}
\def\ACC{{\cal R}(\C)}
\def\ACCa{{\cal R}(\C _1)}
\def\ACCb{{\cal R}(\C _2)}
\def\rACCd{G_{\C}\textrm{-}\mathrm{dim}_R\ }
\def\AsCCd{G_{\C}\textrm{-}\mathrm{dim}_S\ }
\def\ACCd{G_{\C}\textrm{-}\mathrm{dim}\ }

\def\gca{G_{C _1}\textrm{-}\mathrm{dim}\ }
\def\gcb{G_{C _2}\textrm{-}\mathrm{dim}\ }
\def\h{\mathrm{H}}
\def\ord{\mathrm{order}}
\def\s{{s}}
\def\ii{{i}}
\def\amp{{a}}
\def\Z{{\mathbb Z}}
\def\qed{$\Box$}
\def\M{M ^{\bullet }}
\def\N{N ^{\bullet }}
\def\X{X ^{\bullet }}
\def\Y{Y ^{\bullet }}
\def\P{P ^{\bullet }}
\def\Q{Q ^{\bullet }}
\def\F{F ^{\bullet }}
\def\C{C ^{\bullet }}

\def\D{D ^{\bullet }}
\newcommand{\I}{{\rm I}}
\newcommand{\PP}{{\rm P}}

\begin{document}


\title{Homological invariants associated to \\
\vspace{12pt} 
semi-dualizing bimodules}

\author{Tokuji Araya, Ryo Takahashi, and Yuji Yoshino}

\date{\empty}

\maketitle


\begin{abstract}
Cohen-Macaulay dimension for modules over a commutative noetherian local ring has been defined by A. A. Gerko.
That is a homological invariant sharing many properties with projective dimension and Gorenstein dimension.
The main purpose of this paper is to extend the notion of Cohen-Macaulay dimension for modules over commutative noetherian local rings to that for bounded complexes over non-commutative noetherian rings.
\end{abstract}


\section{Introduction}
Cohen-Macaulay dimension for modules over a commutative noetherian local ring has been defined by A. A. Gerko \cite{Gerko}.
That is to be a homological invariant of a module which shares a lot of properties with projective dimension and Gorenstein dimension.
The aim of this paper is to extend this invariant of modules to that of chain complexes, even over non-commutative rings.
We try to pursue it in the most general context possible.

The key role will be played by semi-dualizing bimodules, which we introduce in this paper to generalize semi-dualizing modules in the sense of Christensen \cite{Chris}.
The advantage to consider an $(R, S)$-bimodule structure on a semi-dualizing module $C$ is in the duality theorem.
Actually we shall show that $\Hom _R ( - , C)$ (resp. $\Hom _S (-, C)$) gives a duality between subcategories of $\rmod$ and $\mods$.
We take such an idea from non-commutative ring theory, in particular,  Morita duality and tilting theory.

In \S 2 we present a precise definition of a semi-dualizing bimodule and show several properties.
Associated to a semi-dualizing $(R, S)$-bimodule $C$, of most importance is the notion of the $\gc$-dimension of an $R$-module and the full subcategory $\rAC$ of $\rmod$ consisting of all $R$-modules of finite $\gc$-dimensions.
Under some special conditions the $\gc$-dimension will coincide with the Cohen-Macaulay dimension of a module.

In \S 3 we extend these notions to the derived category, hence to chain  complexes.
We introduce the notion of the trunk module of a complex, and as one of the main results of this paper, we shall show that the $\gc$-dimension of a complex is essentially given by that of its trunk module (Theorem \ref{3.4}).
By virtue of this theorem, we can show that many of the assertions concerning $\gc$-dimensions of modules will hold true for $\gc$-dimensions of complexes.

In \S4 we shall show that a semi-dualizing bimodule, more generally a semi-dualizing complex of bimodules, yields a duality between subcategories of the derived categories.
This second main result (Theorem \ref{3.12}) of this paper gives the advantage from considering the bimodule structure of a semi-dualizing module.

In \S 5 we apply the theory to case that base rings are commutative.
Surprisingly, if the both rings $R$, $S$ are commutative, then we shall see that a semi-dualizing $(R, S)$-bimodule is nothing but a semi-dualizing $R$-module, and actually $R=S$ (Lemma \ref{4.1}).
In this case, we are able to apply the subcategories $\rAC$ and the Gorenstein dimension of the dualizing complex to obtain (in Corollary \ref{4.6}) a new characterization of Gorenstein rings.


\section{$\gc$-dimensions for modules}

{\em Throughout the present paper, we assume that $R$ is a left noetherian ring.
Let $\rmod $ denote the category of finitely generated left $R$-modules.
We also assume that $S$ is a right noetherian ring and $\mods $ denotes the category of finitely generated right $S$-modules. 
When we say simply an $R$-module (resp. an $S$-module), we mean a finitely generated left $R$-module (resp. a finitely generated right $S$-module).}

\vspace{12pt}

In this section, we shall define the notion of  $\gc$-dimension of a module, and study its properties.
For this purpose, we begin with defining a semi-dualizing bimodule.

\begin{defn}\label{2.1}
We call an $(R, S)$-bimodule $C$ a {\it semi-dualizing bimodule} if the following conditions hold.
\begin{itemize}
	\item[(1)]  The right homothety $S$-bimodule morphism $S \to \Hom _R (C,C)$ is a bijection.
	\item[(2)]  The left homothety $R$-bimodule morphism $R \to \Hom _S (C,C)$ is a bijection.
	\item[(3)]  $\Ext ^{i}_R(C,C)=0$ for all $i>0$.
	\item[(4)]  $\Ext ^{i}_S(C,C)=0$ for all $i>0$.
\end{itemize}
\end{defn}

{\em In the rest of this section, $C$ always denotes a semi-dualizing $(R, S)$-bimodule.}

\begin{defn}\label{2.3}
We say that an $R$-module $M$ is $C$-{\it reflexive} if the following conditions are satisfied.
\begin{itemize}
\item[(1)] $\Ext ^{i}_R(M,C)=0$ for all $i>0$.
\item[(2)] $\Ext ^{i}_S(\Hom _R(M,C),C)=0$ for all $i>0$.
\item[(3)] The natural morphism $M \to \Hom_S(\Hom _R(M,C),C)$ is a  bijection.
\end{itemize}
\end{defn}

One can of course consider the same for right $S$-modules by symmetry.

\begin{defn}\label{2.4}
If the following conditions hold for $N \in \mods$, we say that $N$ is $C$-reflexive.
\begin{itemize}
\item[(1)] $\Ext ^{i}_S(N,C)=0$ for all $i>0$.
\item[(2)] $\Ext ^{i}_R(\Hom _S(N,C),C)=0$ for all $i>0$.
\item[(3)] The natural morphism $N \to \Hom_R(\Hom _S(N,C),C)$ is a  bijection.
\end{itemize}
\end{defn}

\begin{ex}\label{2.5}
\begin{itemize}
\item[(1)] Both the ring  $R$ and the semi-dualizing module $C$ are $C$-reflexive $R$-modules.
Similarly, $S$ and $C$ are $C$-reflexive $S$-modules.
\item[(2)] Let $R$ be a finite dimensional algebra.
Then every tilting $R$-module is a semi-dualizing module (cf. \cite[(4.1)]{Ringel}).
\item[(3)] Let $R$ be a left and right noetherian ring. 
Then the ring $R$ itself is a semi-dualizing $(R, R)$-bimodule and the $R$-reflexive modules coincide with the modules whose G-dimension is equal to $0$ (cf. \cite[Proposition 3.8]{AB}).
\item[(4)] Let $R$ be a commutative Cohen-Macaulay local ring with dualizing module $K$. 
Then $K$ is a semi-dualizing module and the $K$-reflexive modules coincide with the maximal Cohen-Macaulay modules (cf.\cite[Theorem 3.3.10]{BH}).
\end{itemize}
\end{ex}

\begin{thm}\label{2.6}
\begin{itemize}
\item[$(1)$] Let $0 \to L_1 \to L_2 \to L_3 \to 0$ be a short exact sequence  either in $\rmod$ or in $\mods$.
Assume that $L_3$ is  $C$-reflexive.
Then, $L_1$ is $C$-reflexive if and only if so is $L_2$.

\item[$(2)$] 
If $L$ is a $C$-reflexive module, then so is any direct summand of $L$.
In particular, any projective module is $C$-reflexive.

\item[$(3)$] The functors $\Hom _R (-,C)$ and $\Hom _S (-,C)$ yield a duality between the full subcategory of $\rmod$ consisting of all $C$-reflexive $R$-modules and the full subcategory of $\mods$ consisting of all $C$-reflexive $S$-modules.
\end{itemize}
\end{thm}

\begin{pf}
(1) Let $0 \to L_1 \to L_2 \to L_3 \to 0$ be a short exact sequence in $\rmod$.
Suppose that $L_3$ is $C$-reflexive.
Applying the functor $\Hom _R(-,C)$ to this sequence, we see that the sequence 
$$
0 \to \Hom _R(L_3,C) \to \Hom _R(L_2,C) \to \Hom _R(L_1,C) \to 0
$$
 is exact, and $\Ext ^i_R(L_2,C) \cong \Ext ^i_R(L_1,C)$ for $i>0$.
Now applying the functor $\Hom _S(-,C)$, we will have an exact sequence 
$$
\begin{CD}
0 @>>> \Hom _S(\Hom _R(L_1,C),C) @>>> \Hom _S(\Hom _R(L_2,C),C)\\
@>>> \Hom _S(\Hom _R(L_3,C),C) @>>> \Ext^1 _S(\Hom _R(L_1,C),C)\\
@>>> \Ext^1 _S(\Hom _R(L_2,C),C) @>>> 0
\end{CD}
$$
and the isomorphisms  $\Ext^i _S(\Hom _R(L_1,C),C) \cong \Ext^i _S(\Hom _R(L_2,C),C)$ for $i>1$.
It is now easy to see from the diagram chasing that  $L_1$ is $C$-reflexive if and only if $L_2$ is $C$-reflexive.

(2) Trivial.

(3) It is straightforward to see that both functors send $C$-reflexive modules to $C$-reflexive modules (over the respective rings), and that the compositions of them are the identity functors (for the respective categories).
\qed
\end{pf}

\begin{ex}
Let $R$ be a laft and right noetherian ring.
We denote ${\cal{G}}_R$ the full subcategory of $\rmod$ consisting of all $R$-reflexive (left) $R$-modules and denote ${\cal{G}}_{R^{op}}$ the full subcategory of $\modr$ consisting of all $R$-reflexive (right) $R$-modules.
In this situation, Theorem \ref{2.6}.(3) says that $\Hom _R(-,R)$ gives a duality between ${\cal{G}}_R$ and ${\cal{G}}_{R^{op}}$.
\end{ex}

\begin{lemma}\label{2.8}
The following conditions are equivalent for $M \in \rmod$ and $n\in\Z$.
\begin{itemize}
\item[$(1)$] 
There exists an exact sequence
$$
0 \to X^{-n} \to X^{-n+1} \to \cdots \to X^0 \to M \to 0
$$
such that each $X^i$ is a $C$-reflexive $R$-module.

\item[$(2)$] 
For any projective resolution
$$
P^{\bullet}\ :\ \cdots \to P^{-m-1} \to P^{-m} \to \cdots \to P^0 \to M \to 0
$$
of $M$ and for any $m \geq n$, we have that $\coker (P^{-m-1} \to P^{-m})$ is a $C$-reflexive $R$-module.

\item[$(3)$] 
For any exact sequence
$$
\cdots \to X^{-m-1} \to X^{-m} \to \cdots \to X^0 \to M \to 0
$$
with each $X^i$ being $C$-reflexive, and for any $m \geq n$, we have that $\coker (X^{-m-1} \to X^{-m})$ is a $C$-reflexive $R$-module.
\end{itemize}
\end{lemma}

\begin{pf}
$(1) \Rightarrow (2)$ : Since $P^{\bullet}$ is a projective resolution of $M$, there is a morphism $\sigma ^{\bullet} : P^{\bullet}\to X^{\bullet}$ of complexes over  $R$ : 
$$
\begin{CD}
\cdots @>>> P^{-n-1} @>>> P^{-n} @>>>P^{-n+1} @>>> \cdots @>>> P^0 \\
@. @V{\sigma ^{-n-1}}VV @V{\sigma ^{-n}}VV @V{\sigma ^{-n+1}}VV @. @V{\sigma ^{0}}VV \\
@. 0 @>>> X^{-n} @>>> X^{-n+1} @>>> \cdots @>>> X^{0}
\end{CD}
$$
Taking the mapping cone of $\sigma ^{\bullet}$, we see that there is an exact sequence
$$
\begin{CD}
\cdots @>>> P^{-m} @>>> \cdots @>>> P^{-n}\\
@>>> P^{-n+1}\oplus X^{-n} @>>> P^{-n+2}\oplus X^{-n+1} @>>> \cdots\\
@>>> P^{0}\oplus X^{-1} @>>> X^{0} @>>> 0.
\end{CD}
$$
It follows from a successive use of Theorem \ref{2.6}.(1) that $\coker (P^{-m-1} \to P^{-m})$ is a $C$-reflexive $R$-module for $m \geq n$.

$(2) \Rightarrow (3)$ : Let $m \geq n$, and set $X = \coker (X^{-m-1} \to X^{-m})$ and $P = \coker (P^{-m-1} \to P^{-m})$.
Since $P^{\bullet}$ is a projective resolution of $M$, there is a chain map   $\sigma ^{\bullet} : P^{\bullet}\to X^{\bullet}$ of complexes over $R$: 
$$
\begin{CD}
0 @>>> P @>>> P^{-m+1} @>>>P^{-m+2} @>>> \cdots @>>> P^0 \\
@. @V{\sigma ^{-m}}VV @V{\sigma ^{-m+1}}VV @V{\sigma ^{-m+2}}VV @. @V{\sigma ^{0}}VV \\
 0 @>>> X @>>> X^{-m+1} @>>> X^{-m+2} @>>> \cdots @>>> X^{0}
\end{CD}
$$
Taking the mapping cone of $\sigma ^{\bullet}$, we see that there is an exact sequence
$$
0 \to P \to P^{-m+1}\oplus X \to P^{-m+2}\oplus X^{-m+1} \to \cdots \to P^{0}\oplus X^{-1} \to X^{0} \to 0.
$$
It then follows from Theorem \ref{2.6}.(1) and \ref{2.6}.(2) that $X$ is a $C$-reflexive $R$-module.

$(3) \Rightarrow (1)$ : Trivial.
\qed
\end{pf}

Imitating the way of defining the G-dimension in \cite[Theorem 3.13]{AB}, we make the following definition.

\begin{defn}\label{2.9}
For $M \in \rmod$, we define the $\gc$-dimension of $M$  by  
$$
\rACd M = \inf
\left\{
\begin{array}{r|l}
  & \textrm{there exists an exact sequence of finite length}\\
 n\ &  0 \to X^{-n} \to X^{-n+1} \to \cdots \to X^0 \to M \to 0,\\
  &\textrm{where each }\ X^i\ \textrm{is  a}\ C\textrm{-reflexive}\ R\textrm{-module}.
\end{array}
\right\}
$$
Here we should note that we adopt the ordinary convention that  $\inf \emptyset = + \infty$.
\end{defn}

\begin{rem}\label{2.10}
First of all we should notice that in the case  $R=S=C$,  the $\gc$-dimension is the same as the G-dimension.

Furthermore, comparing with Theorem \ref{4.3} below, we are able to see that the $\gc$-dimension extends the Cohen-Macaulay dimension over a commutative ring  $R$.
More precisely, suppose that $R$ and $S$ are commutative local rings.
If  there is a semi-dualizing $(R, S)$-bimodule, then  $R$ must be isomorphic to $S$  as we will show later in  Lemma \ref{4.1}.
Thus semi-dualizing bimodules are nothing but semi-dualizing $R$-modules in this case.
One can define the Cohen-Macaulay dimension of an $R$-module  $M$ as 
$$
\cm M  = \inf \{ \rACd M \ | \ C \ \textrm{is a semi-dualizing $R$-module} \}.
$$

Let $C_1$ and $C_2$ be semi-dualizing $R$-modules.
And suppose that an $R$-module $M$  satisfies  $\gca M < \infty$ and $\gcb M < \infty$.
Then we can show that   
$\gca M = \gcb M \ (=\depth R - \depth M)$ (c.f. Lemma \ref{4.2}).
In other words, if the rings $R$ and $S$  are commutative, then the value of the $\gc$-dimension is constant for any choice of semi-dualizing modules $C$ whenever it is finite.
But it follows the next example, if $R$ is non-commutative, this is no longer true.

\begin{ex}\label{2.11}
Let $Q$ be a quiver $e_1 \longrightarrow e_2$, and let $R=kQ$ be the path algebra over an algebraic closed field $k$.
Put $P_1= Re_1$, $P_2= Re_2$, $I_1=\Hom_k (e_1R,k)$, and $I_2=\Hom _k(e_2R,k)$.
Then, it is easy to see that the only indecomposable left $R$-modules are $P_1$, $P_2 \ (\cong I_1 )$, and $I_2$, up to isomorphism. 
Putting  $C_1=P_1\oplus P_2 = R$ and $C_2=I_1 \oplus I_2$, 
we note that  $\End _R(C_1) = \End _R(C_2) = R^{op}$ (here, $R^{op}$ is the opposite ring of $R$), and that $C_1$ and $C_2$ are semi-dualizing $(R, R)$-bimodules.
In this case we have that 
$\gca I_2  (=\g \  I_2)= 1$ and $\gcb I_2 =0$, which take different finite values.
\end{ex}

\end{rem}

\begin{thm}\label{2.12}
If $\rACd M < \infty$ for a module  $M \in \rmod$, then 
$$
\rACd M = \sup \{ \ n \ | \ \Ext ^n_R(M,C) \not= 0 \ \}.
$$
\end{thm}

\begin{pf}
We prove the theorem by induction on $\rACd M$.

Assume first that $\rACd M = 0$.
Then $M$ is a $C$-reflexive module, and hence we have $\sup \{ \ n \ | \ \Ext ^n_R(M,C) \not= 0 \ \}=0$ from the definition.

Assume next that $\rACd M =1$.
Then there exists an exact sequence $0 \to X^{-1} \to X^0 \to M \to 0$ 
 where  $X^0$ and  $X^{-1}$  are  $C$-reflexive $R$-modules.
Then it is clear that  $\Ext _R^i(M, C)=0$  for $i >1$.
We must show that  $\Ext _R^1 (M, C) \not= 0$.
To do this, suppose  $\Ext _R^1 (M, C) =0$.
Then we would have an exact sequence 
\begin{equation}\tag{2.11.1}\label{wildturkey}
0 \to \Hom _R(M,C) \to \Hom _R (X^0 , C) \to \Hom _R(X^{-1}, C) \to 0.
\end{equation} 
Then, writing  the functor  $ \Hom_S (\Hom _R(-, C),C)$  as $F$, 
we get from this the commutative diagram with exact rows 
$$
\begin{CD}
0 @>>> X^{-1} @>>> X^0 @>>> M @>>> 0 \\
 @.      @V{\cong}VV @V{\cong}VV  @VVV \\
0 @>>> F (X^{-1}) @>>>  F (X^0)  @>>> F(M)  @>>> 0,  \\
\end{CD}
$$ 
hence the natural map  $M \to F (M)$  is also an isomorphism.
Furthermore, it also follows from (\ref{wildturkey}) that 
$\Ext ^i _R (\Hom _R(M, C), C) = 0$  for  $i >0$.
Therefore we would have  $\rACd M =0$, a contradiction.
Hence  $\Ext ^1_R (M,C) \not= 0$ as desired.

Finally assume that $\rACd M =m> 1$. 
Then there exists an exact sequence $0 \to X^{-m} \to X^{-m+1} \to \cdots \to X^0 \to M \to 0$  such that each $X^i$ is a $C$-reflexive $R$-module.
Putting  $M'=\coker (X^{-2} \to X^{-1})$, we note that the sequence $0 \to M' \to X^0 \to M \to 0$ is exact and $\rACd M' = m-1>0$.
Therefore, the induction hypothesis implies that $\sup \{ \ n \ | \ \Ext ^n_R(M',C) \not= 0 \ \}=m-1$.
Since $\Ext ^n_R(X^0,C) = 0$ for  $n >0$, 
it follows that $\sup \{ \ n \ | \ \Ext ^n_R(M,C) \not= 0 \ \}=m$ as desired.
\qed
\end{pf}

If $R$ is a left and right noetherian ring and if  $R=S=C$, then the equality  $\mathrm{G}_R\textrm{-}\mathrm{dim}_R M = \g \ M$  holds by definition.
We should remark that if  $R$  is a Gorenstein commutative ring, then any $R$-module $M$  has finite G-dimension and it can be embedded in a short exact sequence of the form  $0 \to F \to X \to M \to 0$ with  $\pd F < \infty$ and $\g \ X =0$.
Such a short exact sequence is called a Cohen-Macaulay approximation of $M$.
For the details, see \cite{AB2}.

We can prove an analogue of this result.
To state our theorem, we need several notations from \cite{AB2}.
Now let  $C$  be a semi-dualizing $(R, S)$-bimodule as before.
We denote  by  ${\cal{G}}_C$  the full subcategory of  $\rmod$  consisting of all  $C$-reflexive $R$-modules, and by  $\rAC$  the full subcategory consisting of $R$-modules of finite $\gc$-dimension.
And  $\add (C)$  denotes the subcategory of all direct summands of direct sums of copies of  $C$.
It is obvious that  $\add (C) \subseteq {\cal{G}}_C$ and that the objects of  $\add(C)$  are injective objects in ${\cal{G}}_C$, because  $\Ext ^i_R(X, C) =0$ for  $X \in {\cal{G}}_C$  and  $i >0$.
The following lemma says that  $C$  is an injective cogenerator of  ${\cal{G}}_C$.

\begin{lemma}\label{2.15}
Suppose an $R$-module $X$  is $C$-reflexive, hence  $X \in {\cal{G}}_C$.
Then there exists an exact sequence $0 \to X \to C^0 \to X^1 \to 0$ 
where  $C^0 \in \add (C)$  and $X^1 \in {\cal{G}}_C$.
In particular, we can resolve $X$  by objects in $\add(C)$  as  
$$
0 \to X \to C^{0} \to C^{1} \to C^{2} \to \cdots, \ (C^i \in \add (C)).
$$
\end{lemma}

\begin{pf}
It follows from Theorem \ref{2.6}.(3) that $Y = \Hom _R(X,C)$ is a $C$-reflexive $S$-module.
Take an exact sequence $0 \to Y' \to S^{\oplus n} \to Y \to 0$  to get the
 the first syzygy $S$-module $Y'$ of $Y$.
Applying the functor $\Hom _S (-,C)$, we obtain an exact sequence $0 \to X \to C^{\oplus n} \to \Hom _S (Y',C) \to 0$.
Since $Y'$ is a $C$-reflexive $S$-module, we see that $\Hom _S (Y',C)$ is a $C$-reflexive $R$-module again.
\qed
\end{pf}

To state the theorem, let us denote
$$
\widehat{\add(C)} = \left\{ 
\begin{array}{r|l}
&\textrm{there is an exact sequence of finite length } \\ 
F \in \rmod &0 \to C^{-n} \to C^{-n+1} \to \cdots \to C^0 \to F \to 0 \\
&\textrm{where each } \ C^i \in \add(C)
\end{array}
\right\}.
$$
Then it is easy to prove the following result in a completely similar way to  the proof of \cite[Theorem 1.1]{AB2}.

\begin{thm}\label{appro}
Let  $M \in \rmod$, and suppose $\rACd M < \infty$, hence  $M \in \rAC$. 
Then there exist short exact sequences
\begin{equation}\tag{2.13.1}\label{ccc}
0 \to F_M \to X_M \to M \to 0
\end{equation}
\begin{equation}\tag{2.13.2}\label{f}
0 \to M \to F^M \to X^M \to 0
\end{equation}
where  $X_M$ and $X^M$ are in  $\gc$, and  $F_M$ and $F^M$ are in 
$\widehat{\add(C)}$. 
\end{thm}

\begin{rem}
Let  $X$  be a $C$-reflexive $R$-module.
Since  $\Ext ^i (X, C) =0$ for $i>0$, it follows that  $\Ext ^i (X, F) = 0$ for $F \in \widehat{\add(C)}$  and $i>0$.
Hence, from  (\ref{ccc}),  we have an exact sequence 
$$
0 \to \Hom _R(X, F_M) \to \Hom _R (X, X_M) \to \Hom _R (X, M) \to 0.
$$
This means that any homomorphism from any $C$-reflexive $R$-module $X$  to  $M$ factors through the map  $X_M \to M$.
In this sense, the exact sequence (\ref{ccc}) gives an approximation of  $M$ by the subcategory  $\gc$. 
\end{rem}

\begin{rem}\label{2.16}
We can of course define  $\AsCd N$  for an $S$-module  $N$  as in the same manner as we define $\rAC$-dimension.
And it is clear by symmetry  that it satisfies that 
$\AsCd N = \sup \{ n | \ \Ext ^n_S(N, C) \not= 0\}$ if the former is finite etc. 
\end{rem}


\section{$\gc$-dimensions for complexes}
Again in  this section, we assume that $R$ (resp. $S$) is a left (resp. right) noetherian ring.
We denote by $\DR $ (resp. $\DS $)  the derived category of  $\rmod$  (resp. $\mods$)  consisting of complexes with bounded finite homologies.

For a complex  $\M$  we always write it as 
$$
\cdots 
\to M^{n-1} \overset{\partial_{M}^{n}}{\longrightarrow } 
M^{n} \overset{\partial_{M}^{n+1}}{\longrightarrow } 
M^{n+1} \overset{\partial_{M}^{n+2}}{\longrightarrow } 
M^{n+2} \to \cdots,
$$
and the shifted complex  $\M [m]$  is the complex with 
 $\M [m]^n = M^{m+n}$ and $\partial_{M[m]}^n=(-1)^m\partial_{M}^{m+n}$.

According to Foxby \cite{Fox2}, 
we define the {\it supremum}, the {\it infimum} and the {\it amplitude} of a complex  $\M$  as follows;

\[
\begin{cases}
\s (\M ) = \sup \{ \ n \ | \ \h ^n(\M ) \not= 0\ \},\\
\ii (\M ) = \inf \{ \ n \ | \ \h ^n(\M ) \not= 0\ \},\\
\amp (\M ) = \s (\M )-\ii (\M ).
\end{cases}
\]

Note that  $\h (\M)=0 \iff \s (\M) = -\infty \iff \ii (\M) = + \infty \iff \amp (\M) = -\infty$.

Suppose in the following that $\h (\M) =0$.
A complex  $\M$  is called bounded if  $\s (\M) < \infty$ and  $\ii (\M) > -\infty$ (hence $\amp (\M) < \infty$).
And  $\DR$  is, by definition, consisting of bounded complexes with finitely generated homology modules. (We remark that for each component $M^n$ of $\M \in \DR$ is not necessary finitely generated.)
Thus, whenever  $\M \in \DR$,  we have 
$$
- \infty <  \ii(\M)  \leq \s (\M) < +\infty.
$$
and $\amp (\M)$  is a non-negative integer.

We remark that the category $\rmod $ can be identified with the full subcategory of $\DR $ consisting of all the complexes $\M\in\DR $ with  $\s (\M )=\ii (\M )=\amp (\M )=0$  or otherwise  $\M\cong 0$.
Through this identification we always think of  $\rmod$  as the full subcategory of  $\DR$.

For a complex  $\P$, if each component $P^i$  is a finitely generated projective module, then we say that  $\P$  is a projective complex.
For any complex $\M \in \DR$, 
we can construct a projective complex $\P$ and a chain map $\P \to \M$ that  yields an isomorphism in $\DR$.
We call such  $\P \to \M$  a semi-projective resolution of  $\M$. 
If  $\h (\M) \not= 0$  and  $s = \s (\M)$ is finite, then 
we can take a semi-projective resolution  $\P$  of  $\M$  in the form;
$$
\cdots \to 
P^{s-2} \overset{\partial_{P}^{s-1}}{\longrightarrow } 
P^{s-1} \overset{\partial_{P}^{s}}{\longrightarrow } 
P^{s} \to 0 \to 0 \to \cdots, \ (\textrm{i.e.} \   P^i = 0 \ \textrm{for} \ i >s).
$$
We call such a semi-projective resolution with this additional property  a {\it   standard}  projective resolution of  $\M$.

For a projective complex $\P (\not\cong 0)$  and an integer $n$, we can consider two kinds of truncated complexes:
\[
\begin{cases}
\tau ^{\leq n} \P = 
(\cdots \to P^{n-2} \overset{\partial_{P}^{n-1}}{\longrightarrow } 
P^{n-1} \overset{\partial_{P}^{n}}{\longrightarrow } P^{n} \to 0 \to 0 \to \cdots)
\\
\tau ^{\geq n} \P = 
(\cdots \to 0 \to 0 \to P^{n} \overset{\partial_{P}^{n+1}}{\longrightarrow } P^{n+1} \overset{\partial_{P}^{n+2}}{\longrightarrow } P^{n+2} \to \cdots )\end{cases}
\]

\begin{defn}[$\omega$-operation]\label{omegaoperation}
Let  $\M \in \DR$, $\h(\M)\neq 0$ and $s=\s (\M)$.
Taking a standard projective resolution $\P$  of $\M$, 
we define the projective complex  $\omega \P$  by
\begin{equation}\tag{3.1.1}\label{defomega}
\omega \P = (\tau ^{\leq s-1}\P )[-1].
\end{equation}
\end{defn}

Note from this definition that 
$\omega \P$  and  $\P [-1]$  share the same components in degree $\leq s$.
We can also see from the definition that there is a triangle of the form 
\begin{equation}\tag{3.1.2}\label{trunkseq}
\omega \P  \to P^{s}[-s] \to \M \to \omega \P[1].
\end{equation}
Therefore, if  $\M$  is a module $M \in \rmod$, then  $\omega \P$  is  isomorphic to a first syzygy module of  $M$. 
Note that  $\omega \P$ is not uniquely determined by  $\M$.
Actually it depends on the choice of a standard projective resolution  $\P$, but is unique up to a projective summand in degree $s$. 
It is easy to prove the following lemma. 

\begin{lemma}\label{i and a} 
Let  $\M \in \DR$, $\h (\M) \neq 0$ and let  $\P$  be a standard projective resolution of  $\M$. 
Now suppose that $\amp (\M) > 0$.
Then, 
\begin{itemize}
\item[$(1)$]
$\ii (\omega \P) = \ii (\M) +1$,
\item[$(2)$]
$0 \leq \amp (\omega \P) < \amp (\M)$.
\end{itemize}
\end{lemma}

\begin{pf}
Let  $s = \s (\M)$.
Since the complexes  $\P$  and  $\omega \P [1]$ share the same components in degree $\leq s -1$, we have that 
$\h ^i(\M)=\h ^i (\P) = \h ^{i+1} (\omega \P)$  for $i \leq s-2$  and that
$\h ^{s-1}(\M) = \h ^{s-1} (\P)$  is embedded into $\h ^{s}(\omega \P)$.
The lemma follows from this observation.
\qed
\end{pf}

It follows from this lemma that applying the $\omega$-operation several times to a given projective complex  $\P$, we will have a complex with amplitude $0$, i.e. a shifted module.

\begin{defn}\label{trunkdef}
Let  $\M$  and  $\P$  be as in the lemma.
Then there is the least integer  $b$  with  $\omega ^b \P$ having amplitude $0$.
Thus there is a module  $T \in \rmod$  such that  $\omega ^b \P \cong T[-c]$  for some $c \in \Z$.
We call such a module $T$ the {\it trunk module} of the complex  $\M$. 
\end{defn}

\begin{rem}\label{trunk}
Let  $\M$  and  $\P$  be as in the lemma.
Set  $i = \ii(\M)$, and we see that the trunk module  $T$  is isomorphic to 
 $\tau ^{\leq i} \P [i]$  in  $\DR$, hence  $T  \cong \coker (P^{i-1} \to P^i)$.
Note that the trunk module  $T$  is unique only in the stable category  ${\underline{\rmod}}$.

Note that the integer  $b$  in Definition \ref{trunkdef} is not necessarily  equal to  $\amp (\M)$. 
For instance, consider the complex  $\M = \P = R[2] \oplus R$.
Then  $\amp (\M) = 2$  and  $T = \omega ^1 \P [-1] = R$. 
\end{rem}

Now we fix a semi-dualizing $(R, S)$-bimodule  $C$.
Associated to it, we can consider the following subcategory of  $\DR$.

\begin{defn}\label{3.1}
For a semi-dualizing $(R, S)$-bimodule  $C$, we denote by $\rAC$ the full subcategory of $\DR$ consisting of all complexes $\M$ that satisfy the following two conditions.
\begin{itemize}
\item[(1)] $\RHom _R(\M,C) \in \DS$.
\item[(2)] The natural morphism $\M \to \RHom_S(\RHom _R(\M,C),C)$ is an isomorphism in $\DR$.
\end{itemize}
\end{defn}

If $R$ is a left and right noetherian ring and if  $R=S=C$, 
then we should note from the papers of Avramov-Foxby \cite[(4.1.7)]{AF} and Yassemi \cite[(2.7)]{Yassemi} that  $\rAR = \{ \ \M \in \DR \ | \ \g \ \M < \infty \}$.

First of all we should notice the following fact.

\begin{lemma}\label{triangle}
Let  $C$  be a semi-dualizing $(R, S)$-bimodule as above.
\begin{itemize}
\item[$(1)$] 
The subcategory $\rAC$ of  $\DR$  is a triangulated subcategory which contains $R$, and is closed under direct summands.
In particular, $\rAC$ contains all projective $R$-modules.
\item[$(2)$]
Let $\P$  be a projective complex in  $\DR$.
Then,   $\P \in \rAC$  if and only if $\omega \P \in \rAC$. 
\item[$(3)$] Let  $\M \in \DR$  and let  $T$  be a trunk module of $\M$.
Then  $\M \in \rAC$  if and only if  $T \in \rAC$.
\end{itemize}
\end{lemma}

\begin{pf}
The proof of (1) is standard, and we omit it. 
For (2) and (3), in the triangle (\ref{trunkseq}), noting that $P[-s] \in \rAC$ and that $\rAC$  is a triangulated category, we see that  
$\P \in  \rAC$  is equivalent to that  $\omega \P \in \rAC$.
Since  $T \cong \omega ^b \P [c]$  as in Definition \ref{trunkdef}, this is also equivalent to that  $T \in \rAC$.
\qed
\end{pf}

The following lemma says that $R$-modules in  $\rAC$  form the subcategory of modules of finite  $\gc$-dimension.

\begin{lemma}\label{3.2}
Let  $M$  be an  $R$-module. 
Then the following two conditions are equivalent.
\begin{itemize}
\item[$(1)$] $\rACd M < \infty$,
\item[$(2)$] $M \in \rAC$.
\end{itemize}
\end{lemma}

\begin{pf}
$(1) \Rightarrow (2)$: 
Note from the definition that every $C$-reflexive module belongs to  $\rAC$.
Since $\rACd M < \infty$, there is a finite exact sequence 
 $0 \to X^{-n} \to X^{-n+1} \to \cdots \to X^0 \to M \to 0$ 
where each $X^i$ is $C$-reflexive.
Since  each $X^i$  belongs to  $\rAC$  and since  $\rAC$  is closed under making triangles, we see that  $M \in \rAC$.

$(2) \Rightarrow (1)$: 
Suppose that  $M \in \rAC$ and let  $\P \in \DR$ be a (standard) projective resolution of  $M$.
Since  $\RHom _R(M, C)$  is a bounded complex, it follows that  $s = \s (\RHom _R(M, C))$ is a (finite) non-negative integer.
Since  the complexes  $\Hom _R(\omega ^s \P , C)$  and 
$\Hom _R(\P [-s], C)$  share the same component in non-negative degree,  
we see that 
$\h ^{i}(\RHom _R(\omega ^s \P , C)) = \h^{i+s}(\RHom _R (\P, C)) = 0$ for $i \geq 1$.
Noting that $\omega^s \P$  is isomorphic to the $s$-th syzygy module $\Omega^s M$ of  $M$,  we see from this that $\Ext ^i (\Omega^s M, C)=0$  for  $i>0$. It follows from above lemma, we have  $\omega ^s \P \in \rAC$ and the natural map 
$\Omega ^s M \to \RHom _S ( \Hom _R(\Omega ^s M, C), C)$  is an isomorphism, equivalently 
$\Omega ^sM \cong \Hom _S( \Hom _R(\Omega ^s M, C), C)$  and 
$\Ext ^i ( \Hom _R(\Omega ^s M, C), C) =0$  for $i>0$.
Consequently, we see that  $\Omega ^s M$ is a $C$-reflexive $R$-module, hence  $\rACd  M \leq s < \infty$. 
\qed
\end{pf}

Recall from Theorem \ref{2.12} that if an $R$-module  $M$  has finite $\gc$-dimension, then we have  $\rACd M = \s (\RHom _R (M, C))$.
Therefore it will be reasonable to make the following definition. 

\begin{defn}
Let  $C$  be a semi-dualizing $(R, S)$-bimodule and let  $\M$  be a complex in  $\DR$.
We define the $\gc$-dimension of  $\M$ to be 
$$
\begin{cases}
\rACd \M = \s (\RHom _R(\M, C))  \quad &\textrm{if} \quad  \M \in \rAC, \\
\rACd \M = +\infty        &\textrm{if}  \quad  \M \not\in \rAC.
\end{cases}
$$
\end{defn}
Note that this definition is compatible with that of $\gc$-dimension for $R$-modules in \S 2.
Just noting an obvious equality  
$$
\s (\RHom _R (\M [m],C)) = \s (\RHom _R (\M,C)) +m
$$
for $\M \in \DR$  and  $m \in \Z$, we have the following lemma.

\begin{lemma}\label{shift}
Let  $\M$  be a complex  in $\DR$ and let  $m$  be an integer.
Then we have
$$
\rACd \M[m] = \rACd \M + m.
$$
\end{lemma}

\begin{lemma}\label{3.5}
Let $\M$  be a complex  in $\DR$. 
Then the following inequality holds:
$$
\rACd \M + \ii (\M) \geq 0.
$$
\end{lemma}

\begin{pf}
If  $\M \cong 0$,  then since  $\ii (\M) = + \infty$, the inequality holds obviously.
We may thus assume that  $\h (\M) =0$.
If $\M \not\in \rAC$, then $\rACd \M = +\infty$ by definition, and there is nothing to prove.
Hence we assume  $\M \in \rAC$.
In particular, we have $\M \cong \RHom _S (\RHom _R(\M,C),C)$.
Therefore we have that 
$$
\begin{array}{rcl}
\ii (\M) & = & \ii (\RHom _S (\RHom _R(\M,C),C)) \\
& \geq & \ii(C) - \s (\RHom _R(\M,C)) \\
& = & -\s (\RHom _R(\M,C)) \\
& = & -\rACd \M.
\end{array}
$$
(For the inequality see Foxby \cite[Lemma 2.1]{Fox1}.)
\qed
\end{pf}

\begin{prop}\label{3.01}
For a given complex  $\M \in \DR$, suppose that  $\amp (\M)>0$.
Taking a standard projective resolution  $\P$  of  $\M$, 
we have an equality
$$
\rACd \M = \rACd \omega \P +1.
$$
\end{prop}

\begin{pf}
Note from Lemma \ref{triangle}(2) that 
 $\rACd \M < \infty$  if and only if  $\rACd \omega \P < \infty$. 
Assume that  $n = \rACd \M = \s (\RHom _R(\P, C)) < \infty$  and let  $s = \s (\M)$.
We should note from Lemma \ref{3.5} that 
$$
\begin{array}{rcl}
n+s & = & \rACd \M + \s (\M) \\
& > & \rACd \M + \ii (\M) \\
& \geq & 0.
\end{array}
$$

Since the complex  $\Hom _R(\omega \P , C)$ shares the components in degree $\geq -s$ with  $\Hom _R(\P, C)[1]$, 
we see that 
$\h ^i (\Hom _R(\omega \P, C)) = \h ^{i+1} (\Hom _R(\P, C))$  for  $i > -s$.
Since $n > -s$ as above, it follows that 
$\s (\Hom _R(\omega \P, C)) = \s (\Hom _R(\P, C))-1$.
\qed
\end{pf}

As we show in the next theorem, the $\gc$-dimension of a complex is essentially the same as that of its trunk module.
In that sense, every argument concerning $\gc$-dimension of complexes  will be reduced to that of modules.

\begin{thm}\label{3.4}
Let   $T$  be the trunk module of a complex $\M \in \DR$  as in Definition \ref{trunkdef}.
Then there is an equality 
$$
\rACd \M = \rACd T - \ii (\M).
$$
\end{thm}

\begin{pf}
If  $\M \not\in \rAC$, then the both sides take infinity and the equality holds.
We assume that  $\M \in \rAC$ hence $\rACd \M < \infty$.

We prove the equality by induction on  $\amp (\M)$.
If  $\amp (\M) = 0$  then  $\M \cong  T [-i]$  for the trunk module  $T$  and for $i = \ii (\M)$.
Therefore it follows from Lemma \ref{shift}  
$\rACd \M = \rACd T -i$. 

Now assume  that  $\amp (\M ) >0$, and let $\P$  be a standard projective resolution of  $\M$.
Noting from Lemma \ref{i and a} that we can apply the induction hypothesis on $\omega \P$, we get the following equalities from the previous proposition.
$$
\begin{array}{rcl}
\rACd \M & = & \rACd \omega \P +  1 \\
& = & \rACd T - \ii (\omega \P) +1 \\
& = & \rACd T - \ii (\P) \\
& = & \rACd T - \ii (\M).
\end{array}
$$
\qed
\end{pf}

As one of the applications of this theorem, we can show the following theorem that generalizes Lemma \ref{2.8} to the category of complexes.

\begin{thm}\label{res}
Let  $\M$  be a complex in  $\DR$.
Then the following conditions are equivalent.
\begin{itemize}
\item[$(1)$]
$\rACd \M < \infty$, 
\item[$(2)$]
There is a bounded complex  $\X$ consisting of $C$-reflexive modules and there is a chain map  $\X \to \M$  that is an isomorphism in $\DR$.
\end{itemize}
\end{thm}

\begin{pf}
$(2) \Rightarrow (1)$:
Note that every $C$-reflexive $R$-module belongs to $\rAC$ and that 
$\rAC$  is closed under making triangles. 
Therefore any complexes  $\X$  of finite length consisting of $C$-reflexive modules are also in  $\rAC$, hence  $\rACd \X < \infty$.

$(1) \Rightarrow (2)$:
Assume that  $\rACd \M < \infty$ hence  $\M \in \rAC$.
We shall prove by induction on $\amp (\M)$  that the second assertion holds.
If  $\amp (\M ) =0$, then there is an $R$-module  $T$  such that 
$\M \cong T[-i]$  where  $i = \ii (\M)$.
Since  $\rACd T < \infty$, there is a complex 
$$
\X = \left[ \ 0 \to X^{-m} \to \cdots \to X^{-2} \to X^{-1} \to X^0 \to 0\ \right] 
$$
with each  $X^i$  being  $C$-reflexive and a quasi-isomorphism $\X \to T$.
Thus the complex  $\X [-s]$  is the desired complex for  $\M$.

Now suppose  $a = \amp (\M) >0$ and take a standard projective resolution  $\P$  of  $\M$.
As in (\ref{trunkseq}), we have chain maps $\varphi : P^s [-s] \to \M$ and
$\psi : \omega \P \to P^s[-s]$  that make the triangle
$$
\omega \P \overset{\psi}\longrightarrow 
P^s[-s] \overset{\varphi}\longrightarrow 
\M \to \omega \P [1].
$$ 
Since $\amp (\omega \P) < \amp (\M)$, it follows from the induction hypothesis that there is a chain map  $\rho : \X \to \omega \P$  that gives an isomorphism in  $\DR$, where  $\X$  is a complex of finite length with each  $X^i$  being $C$-reflexive.
Thus we also have a triangle 
$$
\X \overset{\psi\cdot\rho}\longrightarrow 
P^s[-s] \overset{\varphi}\longrightarrow 
\M \to \X [1].
$$ 
Now take a mapping cone  $\Y$  of  $\psi \cdot \rho$.
Then it is obvious that  $\Y$  has finite length and each modules in  $\Y$  is $C$-reflexive, since  $Y^i$  is a module  $X^i$  with at most directly summing  $P^s$.
Furthermore it follows from the above triangle that there is a chain map  $\Y \to \M$  that yields an isomorphism in  $\DR$.
\qed
\end{pf}

Also in the category $\DS$, we can construct the notion similar to that in $\DR$.
\begin{defn}\label{3.7}
Let  $C$  be a semi-dualizing $(R, S)$-bimodule.
We denote by $\AsC$ the full subcategory of $\DS$ consisting of all complexes $\N$ that satisfy the following two conditions.
\begin{itemize}
\item[(1)] $\RHom _S(\N,C) \in \DR$.
\item[(2)] The natural morphism $\N \to \RHom_R(\RHom _S(\N,C),C)$ is an isomorphism in $\DS$.
\end{itemize}	
\end{defn}

\begin{defn}
Let  $C$  be a semi-dualizing $(R, S)$-bimodule and let  $\N$  be a complex in  $\DS$.
We define the $\AsC$-dimension of  $\N$ to be 
$$
\begin{cases}
\AsCd \N = \s (\RHom _S(\N, C))  \quad &\textrm{if} \quad  \N \in \AsC, \\
\AsCd \N = +\infty        &\textrm{if}  \quad  \N \not\in \AsC.
\end{cases}
$$
\end{defn}
Note that all the properties concerning  $\rAC$ and $\rAC$-dimension hold true for  $\AsC$ and  $\AsC$-dimension by symmetry.

\begin{lemma}\label{dual}
Let  $C$  be a semi-dualizing $(R, S)$-bimodule as above.
Then the functors $\RHom _R (-, C)$ and $\RHom _S (-, C)$ yield a duality between the categories  $\rAC$ and $\AsC$.
\end{lemma}

We postpone the proof of this lemma until Theorem \ref{3.12} in the next section, where we prove the duality in more general setting.
Using this lemma we are able to prove the following theorem, which generalizes Theorem \ref{appro}.
We recall that $\add (C)$ is the additive full subcategory of  $\rmod$  consisting of modules that are isomorphic to direct summands of finite direct sums of copies of  $C$.

\begin{thm}\label{appr}
Let  $\M \in \DR$ and suppose that $\rACd \M < \infty$. 
Then there exists  a  triangle
\begin{equation}\tag{3.17.1}\label{c}
\F_M \to \X_M \to \M \to \F_M [1]
\end{equation}
where $\X_M$ is a shifted $C$-reflexive $R$-module, and  $\F_M$ is a complex that is isomorphic to a complex of finite length consisting of    modules    in  $\add (C)$.
\end{thm}

\begin{pf}
Let  $\N = \RHom _R (\M, C)$  and let  $T$  be a trunk module of  $\N$  in the category  $\DS$.
We have a triangle of the following type:
$$
T[-i] \to \P \to \N \to T[-i+1],
$$
where  $i = \ii (\N)$  and $\P$  is a projective $S$-complex of length  $\amp (\N)$.
Note that  $n = \AsCd T$ is finite as well as  $\AsCd \N < \infty$  by Lemma \ref{dual}.
Take the $n$-th syzygy module of  $T$, and we have a $C$-reflexive $S$-module  $U$ with the triangle
$$
U[-i-n] \to \Q \to \N \to U[-i-n+1],
$$
where  $\Q$  is again a projective  $S$-complex of finite length.
Applying the functor  $\RHom _S(-, C)$, we have a triangle
$$
\RHom _S(U, C)[i+n-1] \to \M \to \RHom _S(\Q, C) \to   
\RHom _S(U, C)[i+n]. 
$$
Note that  $\RHom _S(U, C)$  is isomorphic to a $C$-reflexive  $R$-module and that  $\RHom _S(\Q, C)$  is a complex of finite length, each    component of which is a module in  $\add (C)$. 
\qed
\end{pf}

\section{$\gcc$-dimensions for complexes}

The notion of a semi-dualizing bimodule is naturally extended to that of a semi-dualizing complex of bimodules.
For this purpose, let $\C$  be a complex consisting of  $(R,S)$-bimodules and $(R,S)$-bimodule homomorphisms.
Then for a complex  $\M \in \DR$, take an $R$-projective resolution  $\P$ of  $\M$, and we understand  $\RHom _R(\M, \C)$  as the class of complexes of  $S$-modules that are isomorphic in $\DS$  to the complex  $\Hom _R(\P, \C)$.
In this way,  $\RHom _R(-, \C)$  yields a functor  $\DR \to \DS$.
Likewise, $\RHom _S (-, \C)$  yields a functor $\DS \to \DR$.

Let  $s \in S$.
Then we see that the right multiplication  $\rho (s) : \C \to \C$  is a chain map of $R$-complexes. 
Take a projective resolution  $\P$  of  $\C$  as a complex in  $\DR$
 and  a chain map  $\psi : \P \to \C$  of  $R$-complexes.
Combining these two, we have a chain map  $h (s) = \rho (s) \cdot \psi : \P \to \C$, which defines an element of degree 0 in the complex $\Hom _R(\P, \C)$.
In such a way, we obtain the morphism  $h : S \to \RHom _R(\C, \C)$  in  $\DS$, which we call the right homothety morphism.
Likewise, we have the left homothety morphism 
$R \to \RHom _S(\C, \C)$  in  $\DR$.

\begin{defn}\label{3.9}
Let $\C$  be a complex consisting of  $(R,S)$-bimodules and $(R,S)$-bimodule homomorphisms as above.
We call $\C$ a {\it semi-dualizing complex of bimodules} if the following conditions hold.
\begin{itemize}
\item[(1)] The complex  $\C$  is bounded, that is, there are only a finite number of  $i$  with  $\h ^i(\C) \not= 0$.
\item[(2)] The right homothety morphism $S \to \RHom _R (\C ,\C )$ is an isomorphism in  $\DS$.
\item[(3)] The left homothety morphism $R \to \RHom _S (\C ,\C )$  is an isomorphism in  $\DR$.
\end{itemize}
\end{defn}

\begin{defn}\label{3.10}
We denote by $\rACC$ the full subcategory of $\DR$ consisting of all complexes $\M \in \DR$  that satisfy the following conditions.
\begin{itemize}
\item[(1)]  The complex  $\RHom _R(\M, \C)$  of $S$-modules belongs to  $\DS$.
\item[(2)]  The natural morphism $\M \to \RHom_S(\RHom _R(\M, \C), \C)$  gives an isomorphism in  $\DR$.
\end{itemize}

Similarly we can define  $\AsCC$  as the full subcategory  of $\DS$  
consisting of all complexes $\N$ that satisfy the following conditions.
\begin{itemize}
\item[(1')] The complex  $\RHom _S(\N,\C)$  of  $R$-modules belongs to  $\DR$.
\item[(2')]  The natural morphism $\N \to \RHom_R(\RHom _S(\N,\C),\C)$  gives an isomorphism in $\DS$.
\end{itemize}
\end{defn}

\begin{defn}\label{3.11}
\begin{itemize}
\item[(1)]  For a complex $\M \in \DR$, we define the $\rACC$-dimension of $\M$ as 
$$
\rACCd \M =
\begin{cases}
\s (\RHom _R(\M,\C)) \  &\textrm{if}\  \M \in \rACC, \\
+\infty  \  &\textrm{otherwise}. 
\end{cases}
$$
\item[(2)] Similarly we define the $\AsCC$-dimension of a complex  $\N \in \DS$ as
$$
\AsCCd \N =
\begin{cases}
\s (\RHom _S(\N,\C)) \  &\textrm{if}\  \N \in \rACC, \\
+\infty  \  &\textrm{otherwise}. 
\end{cases}
$$
\end{itemize}
\end{defn}

\begin{thm}\label{3.12}
Let  $\C$  be a semi-dualizing complex of $(R, S)$-bimodules. 
Then the functors $\RHom _R (-,\C)$ and $\RHom _S (-,\C)$ give rise to a duality between $\rACC$ and $\AsCC$.
\end{thm}

\begin{pf}
It is straightforward to see that both functors send complexes with $\gcc$-dimension finite (over the respective rings), and that the compositions of them are the identity functors (for the respective categories).
\qed
\end{pf}


\section{$\gcc$-dimension in the commutative case}

In this final section of the paper, we shall observe several properties of  $\gc$-dimension in the case when $R$ and $S$ are commutative local rings.
We begin with the following lemma.

\begin{lemma}\label{4.1}
Let $R$ and $S$ be commutative noetherian rings.
Suppose that there exists a semi-dualizing $(R,S)$-bimodules $C$.
Then $R$ is isomorphic to $S$.
\end{lemma}

\begin{pf}
Let $\phi : R \to \Hom _R (C,C)=S$ and $\psi : S \to \Hom _S (C,C)=R$ be the homothety morphisms.
Since $R$ and $S$ are commutative, we see that they are well-defined ring homomorphism and that  $\psi \phi$ (resp. $\phi \psi$) is the identity map on $R$ (resp. $S$).
Hence  $R \cong S$ as desired.
\qed
\end{pf}

In view of this lemma, we may assume that $R$ coincides with $S$ for our purpose of this section.
Thus we may call a semi-dualizing $(R, S)$-bimodule simply a semi-dualizing module.
For a semi-dualizing complex $\C$, we simply write $\ACC$  for  $\rACC$.
Note that  $\ACCd \M$ (in this paper) is the same as $\g _{\C} \M$ in \cite{Chris} and $\mathrm{G}_{\C}\textrm{-}\dim \M$ in \cite{Gerko}.

From now on, we assume that $R$ is a commutative noetherian local ring with unique maximal ideal $\m$ and residue class field $k=R/\m$.
It is known that  $\ACCd \M$  satisfies the Auslander-Buchsbaum-type equality as well as $\g _{R} \M$.

\begin{lemma}\label{4.2}{\rm \cite[Theorem\,3.14]{Chris}}
For $\M \in \ACC$, 
$$
\ACCd \M = \depth R - \depth \M + \s (\C ),
$$
where the depth $\depth \M$  of a complex $\M$  is defined to be $\ii (\RHom (k,\M ))$.
\end{lemma}

We are now able to state the main result of this section.

\begin{thm}\label{4.3}
The following conditions are equivalent for a local ring $(R, \m, k)$.
\begin{itemize}
\item[$(1)$]  $R$ is a Cohen-Macaulay local ring that is a homomorphic image of a Gorenstein local ring.
\item[$(2)$] For any finitely generated $R$-module $M$, there exists a semi-dualizing module $C$ such that $\ACd M < \infty $.
\item[$(3)$] There exists a semi-dualizing module $C$ such that $\ACd k < \infty $.
\item[$(4)$] For any $\M \in \DR $ there exists a semi-dualizing module $C$ such that $\ACd \M < \infty $.
\item[$(5)$] There exists a semi-dualizing module $C$ such that $\AC = \DR $.
\item[$(6)$] The dualizing complex $\D $ exists and there exists a semi-dualizing module $C$ such that $\ACd \D < \infty $.
\end{itemize}
\end{thm}

\begin{pf}
The implications $(5) \Rightarrow (4) \Rightarrow (2) \Rightarrow (3)$  are trivial.

$(3) \Rightarrow (1)$: 
Since $\ACd k < \infty $, we have $\Ext ^n_R(k,C)=0$ for $n\gg0$.
Hence we see that the injective dimension of $C$ is finite.
Therefore $R$ is Cohen-Macaulay.
(It is well-known that a commutative local ring which admits a finitely generated module of finite injective dimension is Cohen-Macaulay.
For example, see  \cite{Roberts}.)
Note that 
$$
\begin{array}{rcl}
\depth C & = & -\ACd C + \depth R + \s (C) \\
& = & \depth R \\
& = & \dim R.
\end{array}
$$
That is to say, $C$ is a maximal Cohen-Macaulay module.
Since the isomorphism
$\Ext ^d _R (\Ext ^d _R (k ,C), C) \cong k$, where  $d = \dim R$, 
holds, one can show that $C$ is the dualizing module of $R$.
The existence of the dualizing module of $R$ implies that $R$ is a homomorphic image of a Gorenstein local ring.
(See Reiten \cite[Theorem (3)]{Reiten} or Foxby \cite[Theorem 4.1]{Fox3}.)

$(1) \Rightarrow (6)$:  It follows from the condition (1) that $R$ admits the dualizing module $K_R$.
Note that $K_R$ is a semi-dualizing module and isomorphic to the dualizing complex in $\DR$.
Hence  $\mathrm{G}_{K_R}\textrm{-}\dim K_R = 0 < \infty$.

$(6) \Rightarrow (5)$: 
We may assume that $\ii (\D )=0$.
Then note that $\depth \D = \dim R$.
It follows from Lemma \ref{4.2} that
$$
\begin{array}{rcl}
\ACd \D & = & \depth R - \depth \D + \s (C) \\
& = & \depth R - \dim R \\
& \leq & 0.
\end{array}
$$
On the other hand, from Lemma \ref{3.5} we have that
$$
\ACd \D = \ACd \D + \ii (\D) \geq 0.
$$
Consequently, we have  $\dim R = \depth R$.
Hence $R$ is Cohen-Macaulay. 
And this implies that $\D$ is isomorphic to the dualizing module  $K_R$ of $R$.
It is obvious that  $K_R$  is a semi-dualizing module and  every maximal Cohen-Macaulay module is  $K_R$-reflexive.
As a result, every $R$-module has finite $\mathrm{G}_{K_R}$-dimension, 
hence  ${\cal R}(K_R)$ contains all $R$-modules.
Then it follows from Theorem \ref{3.4}  that 
${\cal R}(K_R)$  contains all complexes in $\DR$, hence 
${\cal R}(K_R) = \DR$.
\qed
\end{pf}

Similarly to the above theorem, we can get a result 
for semi-dualizing complexes. 

\begin{thm}\label{4.4}
The following conditions are equivalent for a local ring $(R, \m, k)$.
\begin{itemize}
\item[$(1)$]  $R$ is a homomorphic image of a Gorenstein local ring.
\item[$(2)$]  For any $M \in \rmod $, there exists a semi-dualizing complex $\C$ such that $\ACCd M < \infty $.
\item[$(3)$]  There exists a semi-dualizing complex $\C$ such that $\ACCd k < \infty $.
\item[$(4)$]  For any $\M \in \DR $, there exists a semi-dualizing complex $\C$ such that $\ACCd M < \infty $.
\item[$(5)$]  There exists a semi-dualizing complex $\C$ such that $\ACC = \DR $.
\item[$(6)$]  The dualizing complex $\D $ exists.
\end{itemize}
\end{thm}

\begin{pf}
It is easy to prove the implications 
$(1) \Rightarrow (6) \Rightarrow (5) \Rightarrow (4) \Rightarrow (2) \Rightarrow (3) \Rightarrow (6)$. 
The remaining implication $(6) \Rightarrow (1)$ that is the most difficult to prove follows from \cite[Theorem 1.2]{K}.
\qed
\end{pf}

As final part of the paper we discuss a kind of uniqueness property of semi-dualizing complexes.

\begin{thm}\label{4.5} 
Let $\C _1$ and $\C _2$ be semi-dualizing complexes.
Suppose that $\C _1 \in \ACCb$ and $\C _2 \in \ACCa$.
Then $\C _1 \cong \C _2 [a]$ for some $a\in\Z$.
In particular, we have $\ACCa = \ACCb $.
\end{thm}

For the proof this theorem we need the notion of Poincare and Bass series of a complex.

\begin{rem} 
Let $(R,\m,k)$ be a commutative noetherian local ring.
For a complex  $\M \in \DR$, consider two kinds of formal Laurent series in the variable  $t$; 
\begin{eqnarray*}
&\PP_{\M}(t) = \sum_{n \in \Z} \dim _k \h ^{-n}(\M \LT k)\cdot  t^n, \\
&\I ^{\M}(t) =\sum_{n \in \Z} \dim _k \h ^n(\RHom (k,\M ))\cdot t^n.
\end{eqnarray*}
These series are called respectively the Poincare series and the Bass series of $\M$.
As it is shown in Foxby \cite[Theorem 4.1(a)]{Fox1}, the following equality holds for $\M ,\N\in\DR$.
\begin{equation}\tag{5.6.1}\label{pb}
\I ^{\RHom (\M,\N)}(t)=\PP _{\M}(t)\cdot \I ^{\N}(t)
\end{equation}
\end{rem}

\noindent
{\sc Proof of \ref{4.5}}.\ \ 
Since $\C _1 \in \ACCb$, we have $\C _1 \cong \RHom (\RHom (\C _1 , \C _2),\C _2)$.
Hence, we have from (\ref{pb}) that 
$$
\I ^{\C_1}(t)=\PP _{\RHom (\C_1, \C_2)}(t)\cdot \I ^{\C_2}(t).
$$
Likewise, it follows from  
$\C _2 \cong \RHom (\RHom (\C _2 , \C _1),\C _1)$ that 
$$
\I ^{\C_2}(t)=\PP _{\RHom (\C_2, \C_1)}(t)\cdot \I ^{\C_1}(t).
$$

Since $\h (\C _1)\neq 0$, one can check that $\I ^{\C_1}(t)\neq 0$.
Therefore we have
$\PP _{\RHom (\C _1 , \C _2)}(t) \cdot \PP _{\RHom (\C _2 , \C _1)}(t) =1$.
Since $\PP _{\RHom (\C _1 , \C _2)}(t)$ and $\PP _{\RHom (\C _2 , \C _1)}(t)$ are formal Laurent series with non-negatiove coefficients, we have
$$
\begin{array}{rcl}
& & \ord (\PP _{\RHom (\C _1 , \C _2)}(t)) + \ord (\PP _{\RHom (\C _2 , \C _1)}(t)) \smallskip\\
& = & \ord (\PP _{\RHom (\C _1 , \C _2)}(t) \cdot \PP _{\RHom (\C _2 , \C _1)}(t)) \bigskip\\
& & \deg (\PP _{\RHom (\C _1 , \C _2)}(t)) + \deg (\PP _{\RHom (\C _2 , \C _1)}(t)) \smallskip\\
& = & \deg (\PP _{\RHom (\C _1 , \C _2)}(t) \cdot \PP _{\RHom (\C _2 , \C _1)}(t))
\end{array}
$$

Therefore we have $\PP _{\RHom (\C _1 , \C _2)}(t)=t^a$ and $\PP _{\RHom (\C _2 , \C _1)}(t) = t^{-a}$ for some integer $a$.
Thus it follows that
$$
\begin{array}{rcl}
\C _1 & \cong & \RHom (\RHom (\C _1 , \C _2), \C _2) \vspace{6pt}\\
 & \cong & \RHom (R[-a], \C _2) \vspace{6pt}\\
 & \cong & \C _2 [a],
 \end{array}
$$
as desired.
\qed

Finally we have an interesting corollary of this theorem.

\begin{cor}\label{4.6}
Suppose that $R$ admits the dualizing complex $\D$.
Then $R$ is a Gorenstein ring if and only if $\g \ \D < \infty $.
\end{cor}

\begin{pf}
If  $R$  is Gorenstein then  $\D \cong R$  thus  $\g \ \D = \g \ R = 0$. 
Conversely, assume  $\g \D < \infty$. 
Then we have  $\D \in {\cal R}(R)$. 
On the other hand, we have  $R \in {\cal R} (\D)$, more generally   ${\cal R}(\D)$  contains all $R$-modules by the definition of dualizing complex.
Hence it follows from the theorem that  $\D \cong R [a]$  for some $a \in \Z$, which means  $R$  is a Gorenstein ring.
\qed
\end{pf}


\bigskip

{\sc Tokuji Araya}

{\sc Graduate School of Natural Science and Technology,

Okayama University, Okayama 700-8530, Japan}

{\it E-mail address} : \texttt{araya@math.okayama-u.ac.jp}

\medskip

{\sc Ryo Takahashi}

{\sc Graduate School of Natural Science and Technology,

Okayama University, Okayama 700-8530, Japan}

{\it E-mail address} : \texttt{takahasi@math.okayama-u.ac.jp}

\medskip

{\sc Yuji Yoshino}

{\sc Faculty of Science, Okayama University,

Okayama 700-8530, Japan}

{\it E-mail address} : \texttt{yoshino@math.okayama-u.ac.jp}

\end{document}